\documentclass[10pt]{amsart}
\usepackage{amsmath,amssymb,amsfonts,amsthm,amsopn}
\usepackage{latexsym,graphicx}
\usepackage{xcolor}
\usepackage{color, colortbl}
\usepackage{pb-diagram}
\usepackage[title]{appendix}
\usepackage{tikz-cd}
\usepackage{tikz}
\usepackage{enumerate}  
\usepackage{mathtools}

\mathtoolsset{showonlyrefs}
\usepackage{multirow}

\newtheorem{theorem}{Theorem}[section]

\newtheorem{proposition}[theorem]{Proposition}
\newtheorem{definition}[theorem]{Definition}

\theoremstyle{definition}
\newtheorem{example}[theorem]{Example}

\newcommand{\beqa}{\begin{eqnarray*}}
\newcommand{\eeqa}{\end{eqnarray*}}

\DeclareMathOperator*{\Sp}{Sp}

\DeclareMathOperator*{\Mp}{Mp}
\DeclareMathOperator*{\Sym}{Sym}

\DeclareMathOperator*{\GL}{GL}

\DeclareMathOperator*{\Spec}{Spec}
\DeclareMathOperator*{\id}{id}

\newcommand{\field}[1]{\mathbb{#1}}
\newcommand{\bR}{\field{R}}        %  real numbers
\newcommand{\bZ}{\field{Z}}        %  whole numbers
\newcommand{\bC}{\field{C}}        %  complex numbers
\def\la{\lambda}

\def\cF{\mathcal{F}}              % Calligraphic Letters
\def\cS{\mathcal{S}}

\def\cA{\mathcal{A}}

\def\rd{\bR^d}

\def\rdd{{\bR^{2d}}}

\def\<{\left<}
\def\>{\right>}

\def\mv1{M_v^1}

\def\mn{(m,n)}
\def\mn'{(m',n')}

\newcommand{\norm}[1]{\lVert#1\rVert}

\def\Ren{\mathbb{R}^d}

\def\Sn2{S_{2}(L^{2}(\Ren))}
\def\S1{S_{1}(L^{2}(\Ren))}
\def\sig00{\sigma_{0,0}}
\def\la{\langle}
\def\ra{\rangle}

\begin{document}

\begin{abstract} 
Time-frequency representations stemmed in 1932 with the introduction of the Wigner distribution. For most of the 20th century, research in this area primarily focused on defining joint quasi-probability distributions for position and momentum in quantum mechanics. Applications to electrical engineering were soon established with the seminal works of Gabor and the researchers at Bell Labs. In 2012, Bai, Li and Cheng used for the first time metaplectic operators, defined in the middle of 20th century by Van Hove, to generalize the Wigner distribution and unify effectively the most used time-frequency representations under a common framework. This work serves as a comprehensive up-to-date survey on time-frequency representations defined by means of metaplectic operators, with particular emphasis on the recent contributions by Cordero and Rodino, who exploited metaplectic operators to their limits to generalize the Wigner distributions. Their idea provides a fruitful framework where properties of time-frequency representations can be explained naturally by the structure of the symplectic group.
\end{abstract}

\title[Metaplectic time-frequency representations]{Metaplectic time-frequency representations}

\author{Gianluca Giacchi}
\address{Universit\`a della Svizzera Italiana, IDSIA, Faculty of Informatics, Via la Santa 1, 6962 Lugano, Switzerland}
\email{gianluca.giacchi@usi.ch}

\thanks{}
%\date{}
\maketitle

\begin{center}
	{\em Dedicated to Luigi Rodino as a tribute to his mathematical career. }
\end{center}

\section{Notation}\label{sec2}
	We denote by $x\xi=x\cdot\xi$ the standard inner product in $\rd$. The notation $\la\cdot,\cdot\ra$ is used either to denote the sesquilinear inner product of $L^2(\rd)$ and its unique extension to a duality pairing (antilinear in the second component) $\cS'(\rd)\times\cS(\rd)$, where $\cS(\rd)$ is the Schwartz class of rapidly decreasing smooth functions. If $f,g$ are functions on $\rd$, $f\asymp g$ means that there exist $A,B>0$ such that $Af(x)\leq g(x)\leq Bf(x)$ holds for every $x\in\rd$. If $f,g\in\cS'(\rd)$ are tempered distributions, $\bar f$ denotes the complex conjugate of $f$, whereas $f\otimes g(x,y):=f(x)g(y)$ is their tensor product. We work with complex valued functions and distributions, and real matrices. The space of $m\times n$ real matrices is denoted by $\bR^{m\times n}$. If $M\in\bR^{m\times n}$, $M^T\in\bR^{n\times m}$ denotes its transpose. If $M$ is a square matrix, $\det(M)$ is its determinant, and $\GL(d,\bR)$ is the group of $d\times d$ invertible matrices. We write $M\in\Sym(d,\bR)$ if $M$ is a $d\times d$ symmetric matrix, i.e., $M^T=M$. $I_d$ denotes the $d\times d$ identity matrix, whereas $0_d$ denotes the $d\times d$ matrix with all zero entries. 

\section{Introduction}
\subsection{The 20th century}
	Since the seminal work by Wigner \cite{Wigner1932}, published in 1932, where the {\em Wigner distribution}
	\begin{equation}\label{Wignerintro}
		Wf(x,\xi)=\int_{\rd}f(x+t/2)\overline{f(x-t/2)}e^{-2\pi i\xi t}dt, \qquad f\in L^2(\rd), \quad x,\xi\in\rd
	\end{equation}
	is defined in the context of quantum mechanics, a broad spectrum of time-frequency representations emerged in the literature. According to \cite{Wigner1932}, the Wigner distribution was found by Wigner himself and Szilard some years before 1932 for another, unspecified, purpose. In 1948 Ville introduced the Wigner distribution in the context of signal analysis \cite{Ville1948}, alongside another time-frequency representation
	\begin{equation}\label{Ambiguityintro}
		Af(x,\xi)=\int_{\rd}f(t+x/2)\overline{f(t-x/2)}e^{-2\pi i\xi t}dt, \qquad f\in L^2(\rd), \quad x,\xi\in\rd,
	\end{equation}
	 named {\em ambiguity function} after the work by Woodward \cite{Woodward}, who applied it to radar analysis. Just one year after the publication of Wigner's work, Kirkwood defined 
	\begin{equation}
		W_0f(x,\xi)=f(x)\overline{\hat f(\xi)}e^{-2\pi i\xi x}, \qquad f\in L^2(\rd), \quad x,\xi\in\rd,
	\end{equation}
	\cite{Kirkhood1933}, that was named after Rihaczek due to his contribution in 1968 \cite{Rihaczek1968}. Parallel to these works, Gabor proposed time-frequency shifts 
	\begin{equation}
		\pi(x,\xi)g(t)=e^{2\pi i\xi t}g(t-x), \qquad x,\xi\in\rd,
	\end{equation}
	to decompose signals into fundamental atoms, where he considered a Gaussian {\em window} $g$, \cite{Gabor}. The term {\em short-time Fourier analysis} emerged naturally over the decades, and the very first appearance of the {\em short-time Fourier transform} (STFT) 
	\begin{equation}\label{STFTintro}
		V_gf(x,\xi)=\la f,\pi(x,\xi)g\ra=\int_{\rd}f(t)\overline{g(t-x)}e^{-2\pi i\xi t}dt,\qquad f,g\in L^2(\rd), \quad x,\xi\in\rd
	\end{equation}
	is lost in the history, even though this phraseology was known to Flanagan and Golden, who used it in 1966 in the context of synthesized human speech, \cite{FlanaganGolden66}. A more experimental evolution occurred for the {\em spectrogram}, stemmed from the empirical works of the Bell Labs \cite{Koenig1946}, summarized in the book by Potter,  Kopp and Green \cite{potter1947visible}, and later developed in \cite{Oppenheim1970spectrograms,rabiner1978digital}. Specifically, the spectrogram is defined as
	\begin{equation}\label{spectrogramIntro}
		{\Spec}_g(f)=|V_gf|^2,\qquad f,g\in L^2(\rd).
	\end{equation}
	A key effort to identify a common thread unifying time-frequency representations is due to Cohen \cite{Cohen1966}, who observed how most of the above-mentioned time-frequency representations can be written as convolutions between the ambiguity function and a distribution, later named {\em kernel} by Claasen and Mecklenbr\"auker \cite{claasen1980time}. The family of such time-frequency representations is known as {\em Cohen's class} at the present. 
	It is straightforward to observe that the Cohen's class can be defined in terms of the Wigner distribution, \cite{claasen1980time},
	so that a time-frequency representation $Q$ belongs to the Cohen's class if 
	\begin{equation}
		Qf=Wf\ast \sigma,\qquad f\in\cS(\rd)
	\end{equation}
	for some $\sigma\in\cS'(\rd)$. Besides the Wigner distribution, the ambiguity function, and the Rihaczek distribution, other elements of the Cohen's class are discussed in \cite{CW1989,Margenau1961,Page1952,ZAM1990}. 
	
		The survey by Cohen \cite{Cohenreview} and his subsequent work \cite{cohen1995time} mark the conclusion of the {\em classical era} of time-frequency distributions. For the purpose of this work, it is worth to mention how a crucial part of Cohen's work was relating properties of representations in the Cohen's class
		%, perhaps interpreted as quantization rules, 
		to the structure of the corresponding kernels.
		
	Other than their different properties, the main aspect distinguishing \eqref{Wignerintro} and \eqref{Ambiguityintro} from \eqref{STFTintro} and \eqref{spectrogramIntro} is the presence of a window function $g$ in the latter. In 1980, Szu and Blodgett \cite{Szu1980} defined the {\em cross-Wigner distribution} 
	\begin{equation}\label{CWintro}
		W(f,g)(x,\xi)=\int_{\rd}f(x+t/2)\overline{g(x-t/2)}e^{-2\pi i\xi t}dt, \qquad f,g\in L^2(\rd), \quad x,\xi\in\rd
	\end{equation}
	and the {\em cross-ambiguity function}
	\begin{equation}\label{CAFintro}
		A(f,g)(x,\xi)=\int_{\rd}f(t+x/2)\overline{g(t-x/2)}e^{-2\pi i\xi t}dt, \qquad f,g\in L^2(\rd), \quad x,\xi\in\rd.
	\end{equation}
	According to the authors, the main interest in introducing the window $g$ was in signal decoding: an unknown signal $f$ could be read correctly by the communicating parties, only with the prior knowledge of the window. In 1985, Janssen \cite{janssen1985bilinear} considered new representations in the Cohen's class known nowadays as the {\em$\tau$-Wigner distributions} 
	\begin{equation}\label{Ctauintro}
		W_\tau f(x,\xi)=\int_{\rd}f(x+\tau t)\overline{f(x-(1-\tau)t)}e^{-2\pi i\xi t}dt, \qquad f,g\in L^2(\rd), \quad x,\xi\in\rd,
	\end{equation}
	where $0\leq \tau\leq 1$. Observe that for $\tau=0$, $W_\tau=W_0$ is the Rihaczek distribution, whereas the classical Wigner distribution corresponds to the choice $\tau=1/2$. The cross-$\tau$-Wigner distributions were considered in \cite{boggiatto2010time}.
	
	\subsection{The metaplectic framework}
	{\em Metaplectic operators} stemmed parallel to time-frequency representations. They were defined in 1951 by Van Hove in his Ph.D. thesis \cite{VanHoveLeon1951}, but they comprise operators that were discovered long before, we mention fractional Fourier transforms \cite{FractionalFTs,FractionalFTs2} as non-straightforward examples. 
	In one sentence, the metaplectic group $\Mp(d,\bR)$ is the double cover of the symplectic group $\Sp(d,\bR)$, see \cite[Section 1.1]{cordero2020time} for the precise definition and the notation. Concretely, it is the group of unitary operators on $L^2(\rd)$ generated by
	\begin{align}
		&i^{-d/2}\cF f(\xi)=i^{-d/2}\int_{\rd}f(x)e^{-2\pi i\xi x}dx,\\
		&i^{m}\mathfrak{T}_Mf(x)=i^{m}|\det(M)|^{1/2}f(Mx), \qquad M\in\GL(d,\bR),\\
		&\mathfrak{p}_Qf(x)=e^{i\pi Qx\cdot x}f(x), \qquad Q\in\bR^{d\times d},\quad Q^T=Q,
	\end{align}
	where $f\in\cS(\rd)$, and $m\in\bZ$ is the argument of $\det(M)$ ({\em Maslov index}). For simplicity, we omit the phase factors $i^{-d/2}$ and $i^{m}$ in the discussion below. The projection $\pi^{Mp}:\hat S\in \Mp(d,\bR)\to S\in\Sp(d,\bR)$ is a group homomorphism with kernel $\{\pm \id_{L^2}\}$. 
	
The interplay between the main time-frequency representations and fractional Fourier transforms was already examined in 2001 by Pei and Ding \cite{Pei2001}. However, the idea of generalizing the Wigner distribution using metaplectic operators is first due to Bai, Li and Cheng, who defined the {\em new Wigner--Ville distribution} in \cite{Bai2012}, by replacing $f$ in \eqref{Wignerintro} with $\hat Sf$, $\hat S\in\Mp(1,\bR)$.

	In 2015, Zhang and Luo generalized the new Wigner--Ville transform further \cite{ZhangNewWigner}, improving its ability of detecting linear frequency modulated signals under low Signal-to-Noise Ratio. Their construction can be formalized in dimension $d\geq1$ as follows: fixed $\hat S,\hat S_1,\hat S_2\in\Mp(d,\bR)$, and $\mathfrak{T}_{M_{1/2}}f(x,t)=f(x+t/2,x-t/2)$, 
	\begin{equation}\label{WS1S2S}
		W_{\hat S,\hat S_1,\hat S_2}f=({\id}_{L^2}\otimes\hat S)\mathfrak{T}_{M_{1/2}}(\hat S_1 f\otimes\overline{\hat S_2 f}), \qquad f,g\in L^2(\rd),
	\end{equation}
	where ${\id}_{L^2}\otimes\hat S$ is the unique metaplectic operator so that $({\id}_{L^2}\otimes\hat S)(f\otimes g)=f\otimes \hat Sg$ for $f,g\in L^2(\rd)$, see \cite[Appendix B]{CORDERO2023154}.

	Parallel to the works by Zhang et al, Bayer, Cordero, Gr\"ochenig and Trapasso defined the (cross-){\em matrix Wigner distributions} in 2019 \cite{Bayer2020,CordTrap}. They observed that a peculiar detail shared by the distributions in \eqref{CWintro}, \eqref{STFTintro}, \eqref{CAFintro} and \eqref{Ctauintro} is that they can be written as
	\begin{equation}\label{MatWintro}
		Q(f,g)=\cF_2\mathfrak{T}_M(f\otimes\bar g), \qquad f,g\in L^2(\rd),
	\end{equation} 
	for suitable $M\in\GL(2d,\bR)$ and
	\begin{equation}\label{FT2}
		\cF_2F(x,\xi)=\int_{\rd}F(x,t)e^{-2\pi i\xi t}dt, \qquad F\in\cS(\rdd), \quad x,\xi\in\rd.
	\end{equation}
	
	Up to the factor $|\det(M)|^{-1/2}$, matrix Wigner distributions are defined as in \eqref{MatWintro}. Observe that both $\cF_2$ and $\mathfrak{T}_M$ used in \eqref{MatWintro} are metaplectic operators acting on $L^2(\rdd)$ and, therefore, matrix Wigner distributions are intrinsically constructed through metaplectic operators. 
	
	\subsection{Metaplectic Wigner distributions} 
	The recent contributions by Cordero and Rodino yielded the definition of the $\cA$-Wigner distributions in 2022 \cite{CORDERO202285}. Their core observation was that many of the most important time-frequency representations defined in the last two centuries can be written letting a metaplectic operator on $L^2(\rdd)$ act on the tensors
	\begin{equation}
		f\otimes\bar g(x,y)=f(x)\overline{g(y)}, \qquad f,g\in L^2(\rd), \quad x,y\in \rd,
	\end{equation}
	as outlined in the following examples.
	\begin{example} 
	\begin{enumerate}[(i)]
	\item The STFT can be written as
		\begin{equation}\label{STFT-2}
			V_gf
			=\widehat{\cA_{st}}(f\otimes\bar g), \qquad f,g\in L^2(\rd).
		\end{equation}
		The projection of $\widehat{\cA_{st}}$ can be computed explicitly as
		\begin{equation}
			\cA_{st}
			=\begin{pmatrix}
				I_d & -I_d & 0_d & 0_d\\
				0_d & 0_d & I_d & I_d\\
				0_d & 0_d & 0_d & -I_d\\
				-I_d & 0_d & 0_d & 0_d
			\end{pmatrix}.
		\end{equation}
	\item Similarly, $W_\tau(f,g)=\widehat{\cA_\tau}(f\otimes\bar g)$, where
		\begin{equation}\label{Atau}
			\cA_\tau
			=\begin{pmatrix}
				(1-\tau)I_d & \tau I_d & 0_d & 0_d\\
				0_d & 0_d & \tau I_d & -(1-\tau)I_d\\
				0_d & 0_d & I_d & I_d\\
				-I_d & I_d & 0_d & 0_d
			\end{pmatrix}.
		\end{equation}
		\item The distributions defined by \eqref{WS1S2S} are metaplectic Wigner distributions. The projection can be computed explicitly using the results in \cite[Appendix B]{CORDERO2023154}: let $S=\pi^{Mp}(\hat S)$, $S_j=\pi^{Mp}(\hat S_j)$ ($j=1,2$) have blocks 
		\begin{equation}
			S=\begin{pmatrix}
				A & B\\
				C & D
			\end{pmatrix}, \quad \text{and} \quad 
			S_j=\begin{pmatrix}
				A_j & B_j\\
				C_j & D_j
			\end{pmatrix},
		\end{equation}
		then $W_{\hat S,\hat S_1,\hat S_2}(f,g)=\widehat{\cA_{\hat S,\hat S_1,\hat S_2}}(f\otimes\bar g)$,
		where $W_{\hat S,\hat S_1,\hat S_2}(f,g)$ is the polarized version of \eqref{WS1S2S}, and
		\begin{equation}\label{projWS1S2S}\begin{split}
			\cA_{\hat S,\hat S_1,\hat S_2}&
			=\begin{pmatrix}
				\frac{1}{2}A_1 & \frac{1}{2}A_2 & \frac{1}{2}B_1 & -\frac{1}{2}B_2\\
				AA_1+\frac{1}{2}BC_1 & -AA_2+\frac{1}{2}BC_2 & AB_1+\frac{1}{2}BD_1 & AB_2-\frac{1}{2}BD_2\\
				C_1 & -C_2 & D_1 & D_2\\
				CA_1+\frac{1}{2}DC_1 & -CA_2+\frac{1}{2}DC_2 & CB_1+\frac{1}{2}DD_1  & CB_2-\frac{1}{2}DD_2
			\end{pmatrix}.
		\end{split}
		\end{equation}
	\item Other time-frequency representations that were not discussed in the introduction can be interpreted similarly. This is the case of $\hbar$ -- STFTs, discussed implicitly in \cite{de2015hamiltonian}. For $\hbar>0$, 
	\begin{equation}
		V^\hbar_gf(x,\xi)=(2\pi\hbar)^{-d/2}e^{2\pi i \frac{x\xi}{4\pi\hbar}}V_gf(x,\frac{\xi}{2\pi\hbar}), \qquad f,g\in L^2(\rd),\quad x,\xi\in\rd
	\end{equation} 
	(the constant $(2\pi\hbar)^{-d/2}$ is added here to retrieve a unitary rescaling). Then, $V^\hbar_gf=\widehat{\cA_\hbar}(f\otimes\bar g)$ and its projection is the $4d\times4d$ symplectic matrix
	\begin{equation}
		\cA_\hbar=\begin{pmatrix}
			I_d & -I_d & 0_d & 0_d\\
			0_d & 0_d & 2\pi\hbar I_d & 2\pi \hbar I_d\\
			0_d & 0_d & I_d/2 & -I_d/2\\
			-\frac{1}{4\pi\hbar}I_d & -\frac{1}{4\pi\hbar}I_d & 0_d & 0_d
		\end{pmatrix}.
	\end{equation}
	\end{enumerate}
	\end{example}

	In view of the preceding examples, the {\em (cross-)$\cA$-Wigner distributions} were defined in \cite{CORDERO202285} as in Definition \ref{defMWDs} below. Later, from \cite{Cordero:2024aa} onwards, the term $\cA$-Wigner distribution is used interchangeably with {\em metaplectic Wigner distributions}. 
	
		\begin{definition}\label{defMWDs}
		For a metaplectic operator $\hat\cA\in\Mp(2d,\bR)$ with projection
	\begin{equation}\label{blockAintro}
		\cA=\begin{pmatrix}
			A_{11} & A_{12} & A_{13} & A_{14}\\
			A_{21} & A_{22} & A_{23} & A_{24}\\
			A_{31} & A_{32} & A_{33} & A_{34}\\
			A_{41} & A_{42} & A_{43} & A_{44}
			\end{pmatrix},
	\end{equation}
	the {metaplectic Wigner distribution}, or $\cA$-Wigner distribution, $W_\cA$ is the time-frequency representation 
	\begin{equation}
		W_\cA(f,g)=\hat\cA(f\otimes\bar g), \qquad f,g\in L^2(\rd).
	\end{equation}
	\end{definition}
	Observe that, since $\cA$ only determines $\hat\cA$ up to a phase factor, writing $W_\cA$ instead of $W_{\hat\cA}$ is an abuse of notation. As a direct consequence of the continuity properties of metaplectic operators, we have the following result.
	\begin{theorem}
		Let $W_\cA$ be a metaplectic Wigner distribution.
		\begin{enumerate}[(i)]
			\item $W_\cA:\cS(\rd)\times\cS(\rd)\to\cS(\rdd)$ is continuous.
			\item $W_\cA:L^2(\rd)\times L^2(\rd)\to L^2(\rdd)$ is bounded and {\em Moyal's identity} holds
			\begin{equation}
				\la W_\cA(f_1,g_1),W_\cA(f_2,g_2)\ra = \la f_1,f_2\ra\overline{\la g_1,g_2\ra}, \qquad f_1,f_2,g_1,g_2\in L^2(\rd).
			\end{equation}
			\item $W_\cA:\cS'(\rd)\times\cS'(\rd)\to\cS'(\rdd)$ is continuous.
		\end{enumerate}
	\end{theorem}
		
		\section{Main aspects}
		\subsection{Cohen's class} Examples of metaplectic Wigner distributions belonging to the Cohen's class are the $\tau$-Wigner distributions. It turns out that the property of belonging to the Cohen's class is strictly related to covariance, as it is proved in \cite{CORDERO202285,CORDERO2023109892}.
		
		\begin{definition}\label{covariance}
			A metaplectic Wigner distribution $W_\cA$ is {\em covariant} if for every $z\in\rdd$,
			\begin{equation}
				W_\cA(\pi(z)f,\pi(z)g)=W_\cA(f,g)(\cdot-z), \qquad f,g\in L^2(\rd).
			\end{equation}
		\end{definition}
		Covariant metaplectic Wigner distributions are characterized in \cite[Proposition 2.10]{CORDERO2023109892}	 in terms of their projection. 
		\begin{theorem}\label{thmCovariance}
			Let $W_\cA$ be any metaplectic Wigner distribution having projection $\cA\in\Sp(2d,\bR)$ with block decomposition \eqref{blockAintro}. The following statements are equivalent.
			\begin{enumerate}[(i)]
			\item $W_\cA$ is covariant.
			\item We have that $A_{13},A_{21}\in\Sym(d,\bR)$ and
			\begin{equation}\label{blockAcov}
				\cA=\begin{pmatrix}
					A_{11} & I_d-A_{11} & A_{13} & A_{13}\\
					A_{21} & -A_{21} & I_d-A_{11}^T & -A_{11}^T\\
					0_d & 0_d & I_d & I_d\\
					-I_d & I_d & 0_d & 0_d
				\end{pmatrix}.
			\end{equation}
			\item $W_\cA(f,g)=W(f,g)\ast \cF^{-1}\Phi_{-B_\cA}$, where
			\begin{equation}\label{defBA}
				B_\cA = \begin{pmatrix}
					A_{13} & I_d/2-A_{11}\\
					I_d/2-A_{11}^T & -A_{21}
				\end{pmatrix}.
			\end{equation}
			\end{enumerate}
		\end{theorem}
		
		\subsection{Shift-invertiblity}	
		Other important characters in time-frequency analysis are {\em modulation spaces}, stemmed by the works of Feichtinger \cite{HGfeichtinger,Fei1981}. We refer to \cite{cordero2020time, grochenig} as exhaustive references for the theory of modulation spaces quickly discussed below. Let $m$ be a $v$-moderate weight on $\rdd$,  $0<p,q\leq\infty$ and $g\in\cS(\rd)\setminus\{0\}$. We set $\norm{f}_{M^{p,q}_m}:=\norm{V_gf}_{L^{p,q}_m}$, where $L^{p,q}_m(\rdd)$ are the mixed-norm Lebesgue spaces, and $V_gf$ is the STFT in \eqref{STFTintro}. Different windows $g$ yield to equivalent quasi-norms. In this section, we discuss which metaplectic Wigner distributions $W_\cA$, other than the STFT, satisfy
		\begin{equation}\label{WAMpq}
		\norm{W_\cA(f,g)}_{L^{p,q}_m}\asymp\norm{f}_{M^{p,q}_m}.
		\end{equation}
		
		For a symplectic matrix $\cA\in\Sp(2d,\bR)$ with blocks \eqref{blockAintro}, we identify the submatrix
		\begin{equation}\label{defEA}
			E_\cA=\begin{pmatrix}
				A_{11} & A_{13}\\
				A_{21} & A_{23}
			\end{pmatrix}.
		\end{equation}
		Let $f\in\cS'(\rd)$ and $g_1,g_2\in\cS(\rd)\setminus\{0\}$. The computations in \cite[Theorem 2.22]{CORDERO2023109892} show that for $1\leq p\leq\infty$ and $v_s=(1+|\cdot|^2)^{s/2}$, $s\geq0$,
		\begin{align}
			\norm{f}_{M^p_{v_s}}\lesssim \norm{|W_\cA(f,g_1)| \ast |W_\cA(g_1,g_2)|(-E_\cA\cdot)}_{L^p_{v_s}}.
		\end{align}
		If $E_\cA$ is invertible, then 
		\begin{align}
			\norm{|W_\cA(f,g_1)| \ast |W_\cA(g_1,g_2)|(-E_\cA\cdot)}_{L^p_{v_s}}\asymp\norm{|W_\cA(f,g_1)| \ast |W_\cA(g_1,g_2)(-\cdot)|}_{L^p_{v_s}},
		\end{align}
		since $v_s\asymp v_s\circ E_\cA$. By Young's inequality,
		\begin{equation}
			\norm{f}_{M^p_{v_s}}\lesssim \norm{W_\cA(f,g_1)}_{L^p_{v_s}}\norm{W_\cA(g_1,g_2)}_{L^1_{v_s}},
		\end{equation}
		and we get the upper bound
		\begin{equation}
			\norm{f}_{M^p_{v_s}}\lesssim \norm{W_\cA(f,g_1)}_{L^p_{v_s}}.
		\end{equation}
		
		In \cite{CORDERO2023109892} the authors conjectured condition $E_\cA\in\GL(2d,\bR)$ to be crucial in the identification of those metaplectic Wigner distributions so that \eqref{WAMpq} holds.
		
		\begin{definition}
			We say that $W_\cA$, or equivalently $\cA$, is {\em shift-invertible} if $E_\cA\in\GL(2d,\bR)$.% in \eqref{blockAintro}.
		\end{definition}
		
		Shift-invertibility was discovered for matrix Wigner distributions in \cite{CordTrap}, but it was formalized only in \cite{CORDERO2023109892} in the broader context of $\cA$-Wigner distributions.
		
		\begin{example}
			\begin{enumerate}[(i)]
				\item The Wigner distribution and the STFT are shift-invertible. Moreover, the $\tau$-Wigner distributions are shift-invertible, for every $\tau\neq0,1$.
				\item Shift-invertible matrix Wigner distributions were studied in \cite{Cordero:2024aa}.
				\item A distribution $W_{\hat S,\hat S_1,\hat S_2}$ in \eqref{WS1S2S} is shift-invertible if and only if the projection $S$ of $\hat S$ is free. Indeed, a straightforward computation shows that
				\begin{equation}
				E_{\cA_{\hat S,\hat S_1,\hat S_2}}=\begin{pmatrix} \frac{1}{2}A_1 & \frac{1}{2}B_1 \\ AA_1+\frac{1}{2}BC_1 & AB_1+\frac{1}{2}BD_1\end{pmatrix}=\frac{1}{2}\begin{pmatrix} I_d & 0_d\\ -2A & B\end{pmatrix}S_1,
				\end{equation}
				which is invertible if and only if $B\in\GL(d,\bR)$. 
			\end{enumerate}
		\end{example}
		
		In \cite{CORDERO2023154}, the authors proved that shift-invertibility was sufficient (minor requirements aside) for \eqref{WAMpq} to hold for metaplectic Wigner distributions in the Banach setting, by direct estimating $\norm{W_\cA(f,g)}_{L^{p,q}_m}$ from above and below. Later, in  \cite{CORDERO2024101594}, the same authors retrieved the characterization of shift-invertible distributions as rescaled STFTs. To be exhaustive, if $\cA$ has blocks \eqref{blockAintro}, we need to define the symmetric matrix 
		\begin{equation}
			M_\cA=\begin{pmatrix}
				A_{11}^TA_{31}+A_{21}^TA_{41} & A_{11}^TA_{33}+A_{21}^TA_{43}\\
				A_{13}^TA_{31}+A_{23}^TA_{41}+I_d & A_{13}^TA_{33}+A_{23}^TA_{43}
			\end{pmatrix},
		\end{equation}
		and the symplectic matrix (see \cite{CORDERO2024101594})
		\begin{equation}
			G_\cA=\begin{pmatrix}
			0_d & I_d\\
			I_d & 0_d
		\end{pmatrix}E_\cA^{-1}\begin{pmatrix}
			A_{12} & A_{14}\\
			A_{22} & A_{24}
		\end{pmatrix}.
		\end{equation}
		
		\begin{proposition}
			Let $W_\cA$ be a metaplectic Wigner distribution, and $\cA$ be the corresponding projection having blocks \eqref{blockAintro}. The following statements are equivalent.
			\begin{enumerate}[(i)]
			\item $W_\cA$ is shift-invertible.
			\item For every $f,g\in L^2(\rd)$,
			\begin{equation}\label{Washiftinv}
				W_\cA(f,g)(z)=c_\cA|\det(E_\cA)|^{-1/2}\Phi_{M_\cA}(E_\cA^{-1}z)V_{\widehat{\delta_\cA}g}f(E_\cA^{-1}z), \qquad z\in\rdd,
			\end{equation}
			for some constant $c_\cA\in\bC$, $|c_\cA|=1$, where $\widehat{\delta_\cA}=\cF \widehat{\overline{G_\cA}}$, where ${\widehat{\overline{G_\cA}}}f=\widehat{{G_\cA}}\bar f$.
			\item There exist $E\in\GL(2d,\bR)$, $C\in\Sym(2d,\bR)$, $\hat \delta\in\Mp(d,\bR)$ and $c\in\bC$, $|c|=1$, such that
			\begin{equation}
				W_\cA(f,g)(z)=c|\det(E)|^{1/2}\Phi_{C}(Ez)V_{\widehat{\delta}g}f(Ez), \qquad f,g\in L^2(\rd),\quad z\in\rdd.
			\end{equation}
			\end{enumerate}
		\end{proposition}
		
			The set of $4d\times4d$ shift-invertible symplectic matrices is studied in \cite{Cordero:2024ab} and \cite{Giacchi2024aa}. 
			
		\begin{theorem}\label{shinve-mod}
			Let $W_\cA$ be a metaplectic Wigner distribution, $\cA$ be the corresponding projection having blocks \eqref{blockAintro}, and $E_\cA$ be the submatrix of $\cA$ defined as in \eqref{defEA}. Let the $v$-moderate weight $m$ satisfy $m\circ E_\cA\asymp m$.
			\begin{enumerate}[(i)]
				\item The equivalence \eqref{WAMpq} holds for every $0<p=q\leq\infty$.
				\item If $E_\cA$ is upper block triangular ($A_{21}=0_d$), then \eqref{WAMpq} holds for every $0<p,q\leq\infty$.
			\end{enumerate}
			Moreover, if $p,q\geq1$, the window $g$ in \eqref{WAMpq} can be chosen in $M^1_v(\rd)\setminus\{0\}$.
		\end{theorem}
		
		A result similar to Theorem \ref{shinve-mod} can be formulated for Wiener amalgam spaces $W(\cF L^p_{m_1},L^q_{m_2})(\rd)=\cF(M^{p,q}_{m_1\otimes m_2}(\rd))$, see \cite[Corollary 3.12]{CORDERO2023154} and \cite[Corollary 7.2]{CORDERO2024101594}. These spaces were defined by Feichtinger in 1981 \cite{feichtinger1981banach}.\newline
		
		Theorem \ref{shinve-mod} tells that shift-invertibility is sufficient, along with other minor conditions, to characterize modulation spaces. On the contrary, the Rihaczek distribution is not shift-invertible with 
		\begin{equation}
			\norm{W_0(f,g)}_{L^{p,q}}=\norm{f}_p\norm{\hat g}_q,
		\end{equation}
		and right-hand side is not proportional to $\norm{f}_{M^{p,q}}$, unless $p=q=2$. This fact is more general. In \cite{Giacchi2024aa} it is proved that $M^p$ cannot be characterized by means of non shift-invertible distributions.
		
		\begin{theorem}
		Let $W_\cA$ be a non shift-invertible metaplectic Wigner distributions, and $g\in\cS(\rd)\setminus\{0\}$. Then, if $0<p\leq\infty$, $p\neq2$,
		\begin{equation}
			\{f\in\cS'(\rd):\norm{W_\cA(f,g)}_{{p}}<\infty\}\neq M^{p}(\rd).
		\end{equation}
		\end{theorem}
		
		\subsection{Generalized spectrograms}
		Generalized spectrograms were defined by Boggiatto, De Donno and Oliaro in 2007 \cite{Boggiatto2007}: for fixed $\varphi,\psi\in\cS'(\rd)$, and $f,g\in\cS'(\rd)$,
		\begin{equation}
			{\Spec}^{\varphi,\psi}(f,g)=V_\varphi f\overline{V_\psi g}.
		\end{equation}
		These distributions are metaplectic Wigner distributions under appropriate choices of $\varphi$ and $\psi$. The reader may find a complete characterization of metaplectic Wigner distributions that define generalized spectrograms in \cite[Theorem 4.6]{Cordero:2025aa}.
		
		\subsection{Uncertainty principles}
Uncertainty principles for metaplectic Wigner distributions have been studied independently by many authors. Most of the known results are restricted to particular instances of these time-frequency representations, where the computations are explicit, see the early survey \cite{Giacchi2022metaplectic} and \cite{WangMatrix}. In \cite{grochenig2024benedicks}, the authors formulate Benedick's uncertainty principle for metaplectic Wigner distributions, characterizing which time-frequency representations satisfy 
\begin{center}
	{\em $W_\cA (f,g)$ supported on a set of finite measure $\Rightarrow$ $f=0$ or $g=0$}
\end{center}
in terms of the Iwasawa decomposition of the corresponding projection $\cA$.
		
		\section{Inversion formula and Gabor frames}
		The content of this section is contained in \cite{CORDERO2024101594}. We refer the reader to \cite{cordero2020time,grochenig} for the theory of (Gabor) frames. For every metaplectic Wigner distribution $W_\cA$ there exists a family of continuous operators $\pi_\cA(z):\cS(\rd)\to\cS'(\rd)$, $z\in\rdd$, so that
		\begin{equation}
			W_\cA(f,g)(z)=\la f,\pi_\cA(z)g\ra, \qquad f,g\in \cS(\rd)
		\end{equation}
		(\emph{metaplectic atoms}).
		\begin{theorem}
			If $f\in L^2(\rd)$ and $g,\gamma\in\cS(\rd)$ are so that $\la g,\gamma\ra\neq0$, then 
			\begin{equation}
				f\overset{L^2}{=}\frac{1}{\la\gamma,g\ra}\int_{\rdd}W_\cA(f,g)(z)\pi_\cA(z)\gamma dz.
			\end{equation}
		\end{theorem}
		In general, metaplectic atoms map $\cS(\rd)$ to $\cS'(\rd)$. However, metaplectic atoms of shift-invertible metaplectic Wigner distributions define surjective quasi-isometries of $L^2(\rd)$. In this case, we can consider a system $\mathcal{G}_\cA(g,\Lambda)=\{\pi_\cA(\lambda)g\}_{\lambda\in\Lambda}$, where $g\in L^2(\rd)\setminus\{0\}$ and $\Lambda\subset\rdd$ is discrete. We say that $\mathcal{G}_\cA(g,\Lambda)$ is a {\em metaplectic Gabor frame} (for $L^2(\rd)$) if it is a frame for $L^2(\rd)$. In view of \eqref{Washiftinv}, it is straightforward to relate metaplectic Gabor frames to classical Gabor frames, see \cite[Theorem 6.4]{CORDERO2024101594}. In the case of shift-invertible distributions, the {\em deformation operator} $\widehat{\delta_\cA}$ in \eqref{Washiftinv} plays a central role in the frame theory associated to $W_\cA$. The {\em frame operator} associated to $\mathcal{G}_\cA(g,\Lambda)$ is $S_\cA f=\sum_{\lambda\in\Lambda}W_\cA(f,g)(\lambda)\pi_\cA(\lambda)g$, $f\in L^2(\rd)$, and the {\em canonical dual window} of $g$ is $\gamma_\cA=\widehat{\delta_\cA}^{-1}S_\cA^{-1}\widehat{\delta_\cA}g$. We have that if $m\asymp m\circ E_\cA$, then for every $0<p,q\leq\infty$,
		\begin{equation}
			f\overset{M^{p,q}_m}{=}\sum_{\lambda\in\Lambda}W_\cA(f,g)(\lambda)\pi_\cA(\lambda)\gamma_\cA,
		\end{equation}
		with unconditional convergence if $\max\{p,q\}\neq\infty$, and weak-$\ast$ convergence in $M^{\infty}_{1/v}$, otherwise. Moreover,
		\begin{equation}
			\norm{f}_{M^{p,q}_m}\asymp \norm{(W_\cA(f,g)(\lambda))_{\lambda\in\Lambda}}_{\ell^{p,q}_m},
		\end{equation}
		where $\ell^{p,q}_m(\Lambda)$ are the mixed-norm sequence spaces. We refer the reader to \cite[Theorem 7.3]{CORDERO2024101594} for the precise statement.
		
\bibliographystyle{abbrv}

\begin{thebibliography}{10}

\bibitem{Bai2012}
R.-F. Bai, B.-Z. Li, and Q.-Y. Cheng.
\newblock Wigner-{V}ille distribution associated with the linear canonical
  transform.
\newblock {\em J. Appl. Math.}, 2012(1):740161, 2012.

\bibitem{Bayer2020}
D.~Bayer, E.~Cordero, K.~Gr{\"o}chenig, and S.~I. Trapasso.
\newblock {\em Linear Perturbations of the Wigner Transform and the Weyl
  Quantization}, pages 79--120.
\newblock Springer International Publishing, Cham, 2020.

\bibitem{boggiatto2010time}
P.~Boggiatto, G.~De~Donno, and A.~Oliaro.
\newblock Time-frequency representations of {W}igner type and
  pseudo-differential operators.
\newblock {\em Trans. Am. Math. Soc.}, 362(9):4955--4981, 2010.

\bibitem{Boggiatto2007}
P.~Boggiatto, G.~D. Donno, and A.~Oliaro.
\newblock A class of quadratic time-frequency representations based on the
  short-time {F}ourier transform.
\newblock In J.~Toft, editor, {\em Modern Trends in Pseudo-Differential
  Operators}, pages 235--249, Basel, 2007. Birkh\"auser Basel.

\bibitem{CW1989}
H.-I. Choi and W.~Williams.
\newblock Improved time-frequency representation of multicomponent signals
  using exponential kernels.
\newblock {\em IEEE Trans. Acoust. Speech Signal Process.}, 37(6):862--871,
  1989.

\bibitem{claasen1980time}
T.~Claasen and W.~Mecklenbr\"auker.
\newblock The {W}igner distribution -- a tool for time-frequency signal
  analysis; part {I}{I}{I}: relations with other time-frequency signal
  transformations.
\newblock {\em Philips J. Res.}, 35(6):372--389, 1980.

\bibitem{Cohen1966}
L.~Cohen.
\newblock Generalized phase-space distribution functions.
\newblock {\em J. Math. Phys.}, 7(5):781--786, 05 1966.

\bibitem{Cohenreview}
L.~Cohen.
\newblock Time-frequency distributions-a review.
\newblock {\em Proc. IEEE.}, 77(7):941--981, 1989.

\bibitem{cohen1995time}
L.~Cohen.
\newblock {\em Time-frequency analysis}, volume 778.
\newblock Prentice Hall PTR New Jersey, 1995.

\bibitem{FractionalFTs}
E.~U. Condon.
\newblock Immersion of the {F}ourier transform in a continuous group of
  functional transformations.
\newblock {\em Proc. Natl. Acad. Sci. {USA}}, 23(3):158--164, Mar 1937.

\bibitem{CORDERO2023154}
E.~Cordero and G.~Giacchi.
\newblock Symplectic analysis of time-frequency spaces.
\newblock {\em J. Math. Pures. Appl.}, 177:154--177, 2023.

\bibitem{Cordero:2024ab}
E.~Cordero and G.~Giacchi.
\newblock Excursus on modulation spaces via metaplectic operators and related
  time-frequency representations.
\newblock {\em Sampl. Theory Signal Process. Data Anal.}, 22(1):9, 2024.

\bibitem{CORDERO2024101594}
E.~Cordero and G.~Giacchi.
\newblock Metaplectic {G}abor frames and symplectic analysis of time-frequency
  spaces.
\newblock {\em Appl. Comput. Harmon. Anal.}, 68:101594, 2024.

\bibitem{Cordero:2024aa}
E.~Cordero, G.~Giacchi, and L.~Rodino.
\newblock {Wigner Analysis of Operators. Part {I}{I}: {S}chr\"odinger
  Equations}.
\newblock {\em Comm. Math. Phys.}, 405(7):156, 2024.

\bibitem{Cordero:2025aa}
E.~Cordero, G.~Giacchi, and L.~Rodino.
\newblock A unified approach to time--frequency representations and generalized
  spectrograms.
\newblock {\em J. Fourier Anal. Appl.}, 31(1):9, 2025.

\bibitem{cordero2020time}
E.~Cordero and L.~Rodino.
\newblock {\em Time-Frequency Analysis of Operators}.
\newblock De Gruyter, Berlin, Boston, 2020.

\bibitem{CORDERO202285}
E.~Cordero and L.~Rodino.
\newblock {Wigner analysis of operators. Part {I}: pseudodifferential operators
  and wave fronts}.
\newblock {\em Appl. Comput. Harmon. Anal.}, 58:85--123, 2022.

\bibitem{CORDERO2023109892}
E.~Cordero and L.~Rodino.
\newblock Characterization of modulation spaces by symplectic representations
  and applications to {S}chr\"odinger equations.
\newblock {\em J. Funct. Anal.}, 284(9):109892, 2023.

\bibitem{CordTrap}
E.~Cordero and S.~I. Trapasso.
\newblock Linear perturbations of the {W}igner distribution and the {C}ohen
  class.
\newblock {\em Anal. Appl.}, 18(03):385--422, 2020.

\bibitem{de2015hamiltonian}
M.~A. {de Gosson}.
\newblock Hamiltonian deformations of {G}abor frames: first steps.
\newblock {\em Appl. Comput. Harmon. Anal.}, 38(2):196--221, 2015.

\bibitem{HGfeichtinger}
H.~Feichtinger.
\newblock Modulation spaces on locally compact abelian groups.
\newblock Technical report, Universit{\"a}t Wien, Mathematisches Institut,
  1983.

\bibitem{feichtinger1981banach}
H.~G. Feichtinger.
\newblock {Banach spaces of distributions of Wiener's type and interpolation}.
\newblock In {\em Functional Analysis and Approximation: Proceedings of the
  Conference held at the Mathematical Research Institute at Oberwolfach, Black
  Forest, August 9--16, 1980}, pages 153--165. Springer, 1981.

\bibitem{Fei1981}
H.~G. Feichtinger.
\newblock On a new {S}egal algebra.
\newblock {\em Monatsh. Math.}, 92(4):269--289, 1981.

\bibitem{FlanaganGolden66}
J.~L. Flanagan and R.~M. Golden.
\newblock Phase vocoder.
\newblock {\em Bell System Technical Journal}, 45(9):1493--1509, 1966.

\bibitem{Gabor}
D.~Gabor.
\newblock Theory of communication.
\newblock {\em J. Inst. Electr. Eng. Part III Radio Commun. Eng.},
  93:429p--457, 1946.

\bibitem{Giacchi2024aa}
G.~Giacchi.
\newblock {Boundedness of Metaplectic Operators Within $L^p$ Spaces,
  Applications to Pseudodifferential Calculus, and Time--Frequency
  Representations}.
\newblock {\em J. Fourier Anal. Appl.}, 30(6):69, 2024.

\bibitem{Giacchi2022metaplectic}
G.~Giacchi.
\newblock Metaplectic {W}igner distributions.
\newblock In A.~Tabacco, P.~M., and A.~N., editors, {\em New Trends in Complex
  and Fourier Analysis, and Operator Theory}. Springer Nature, Singapore Pte
  Ltd, to appear.

\bibitem{grochenig}
K.~Gr{\"o}chenig.
\newblock {\em {Foundations of Time-Frequency Analysis}}.
\newblock Applied and Numerical Harmonic Analysis. Birkh{\"a}user, Boston, MA,
  2001.

\bibitem{grochenig2024benedicks}
K.~Gr{\"o}chenig and I.~Shafkulovska.
\newblock {Benedicks-type uncertainty principle for metaplectic time-frequency
  representations}.
\newblock {\em J. Anal. Math.}, 2025.

\bibitem{janssen1985bilinear}
A.~J. E.~M. Janssen.
\newblock Bilinear phase‐plane distribution functions and positivity.
\newblock {\em J. Math. Phys.}, 26(8):1986--1994, 08 1985.

\bibitem{Kirkhood1933}
J.~G. Kirkwood.
\newblock Quantum statistics of almost classical assemblies.
\newblock {\em Phys. Rev.}, 44:31--37, Jul 1933.

\bibitem{Koenig1946}
W.~Koenig, H.~K. Dunn, and L.~Y. Lacy.
\newblock The sound spectrograph.
\newblock {\em J. Acoust. Soc. Am.}, 18(1):19--49, 07 1946.

\bibitem{Margenau1961}
H.~Margenau and R.~N. Hill.
\newblock Correlation between measurements in quantum theory:.
\newblock {\em Prog. Theor. Phys.}, 26(5):722--738, 11 1961.

\bibitem{FractionalFTs2}
V.~Namias.
\newblock The fractional order {F}ourier transform and its application to
  quantum mechanics.
\newblock {\em IMA J. Appl. Math.}, 25(3):241--265, 03 1980.

\bibitem{Oppenheim1970spectrograms}
A.~V. Oppenheim.
\newblock Speech spectrograms using the fast {F}ourier transform.
\newblock {\em IEEE Spectr.}, 7(8):57--62, 1970.

\bibitem{Page1952}
C.~H. Page.
\newblock Instantaneous power spectra.
\newblock {\em J. Appl. Phys.}, 23(1):103--106, 01 1952.

\bibitem{Pei2001}
S.-C. Pei and J.-J. Ding.
\newblock Relations between fractional operations and time-frequency
  distributions, and their applications.
\newblock {\em IEEE Trans Signal Process.}, 49(8):1638--1655, 2001.

\bibitem{potter1947visible}
R.~Potter, G.~Kopp, and G.~H.C.
\newblock {\em Visible Speech}.
\newblock The Bell Telephone Laboratories series. D. Van Nostrand Company,
  1947.

\bibitem{rabiner1978digital}
L.~Rabiner and R.~Schafer.
\newblock {\em Digital Processing of Speech Signals}.
\newblock Prentice-Hall signal processing series. Prentice-Hall, 1978.

\bibitem{Rihaczek1968}
A.~Rihaczek.
\newblock Signal energy distribution in time and frequency.
\newblock {\em IEEE Trans. Inf. Theor.}, 14(3):369--374, 1968.

\bibitem{Szu1980}
H.~H. Szu and J.~A. Blodgett.
\newblock Wigner distribution and ambiguity function.
\newblock {\em AIP Conference Proceedings}, 65(1):355--381, 12 1980.

\bibitem{VanHoveLeon1951}
L.~Van~Hove.
\newblock {\em Sur certaines repr\'esentations unitaires d'un groupe infini de
  transformations}.
\newblock PhD thesis, Acad\'emie royale de Belgique, 1951.

\bibitem{Ville1948}
J.~Ville.
\newblock Th\'eorie et applications de la notion de signal analytique.
\newblock {\em C\^ables Transm.}, 2:61--74, 1948.

\bibitem{WangMatrix}
L.~Wang, M.~Cui, Z.~Qin, Z.~Zhang, and J.~Zhang.
\newblock Matrix-{W}igner distribution.
\newblock {\em Fractal Fract.}, 8(6), 2024.

\bibitem{Wigner1932}
E.~Wigner.
\newblock On the quantum correction for thermodynamic equilibrium.
\newblock {\em Phys. Rev.}, 40:749--759, Jun 1932.

\bibitem{Woodward}
P.~M. Woodward.
\newblock {\em Probability and Information Theory with Applications to Radar}.
\newblock International Series of Monographs on Electronics and
  Instrumentation. Pergamon, second edition, 1953.

\bibitem{ZhangNewWigner}
Z.~Zhang and M.~Luo.
\newblock New integral transforms for generalizing the {W}igner distribution
  and ambiguity function.
\newblock {\em IEEE Signal Process. Lett.}, 22(4):460--464, 2015.

\bibitem{ZAM1990}
Y.~Zhao, L.~Atlas, and R.~Marks.
\newblock The use of cone-shaped kernels for generalized time-frequency
  representations of nonstationary signals.
\newblock {\em IEEE Trans. Acoust. Speech Signal Process.}, 38(7):1084--1091,
  1990.

\end{thebibliography}

\end{document}